\documentclass[a4paper,12pt]{article}

\usepackage{graphicx}      
\usepackage{epsfig} 
\usepackage{amsmath} 
\usepackage{amssymb}  

\newcommand{\D}{\mathcal{D}}

\newcommand{\nocode}[1]{}
  
\title{\bf H$_2$ for HIFOO}

\begin{document}

\author{Denis Arzelier$^1$, Georgia Deaconu$^2$, Suat Gumussoy$^3$, Didier Henrion$^4$}

\footnotetext[1]{CNRS; LAAS; 7 avenue du colonel Roche, F-31077 Toulouse, France;
Universit\'e de Toulouse; UPS, INSA, INP, ISAE; LAAS; F-31077 Toulouse, France (e-mail: denis.arzelier@laas.fr)}
\footnotetext[2]{CNRS; LAAS; 7 avenue du colonel Roche, F-31077 Toulouse, France;
Universit\'e de Toulouse; UPS, INSA, INP, ISAE; LAAS; F-31077 Toulouse, France (e-mail: georgia.deaconu@laas.fr)}
\footnotetext[3]{Katholieke Universiteit Leuven, Department of Computer Science, Belgium
(e-mail: suat.gumussoy@cs.kuleuven.be)}
\footnotetext[4]{CNRS; LAAS; 7 avenue du colonel Roche, F-31077 Toulouse, France;
Universit\'e de Toulouse; UPS, INSA, INP, ISAE; LAAS; F-31077 Toulouse, France. Faculty of Electrical Engineering, Czech Technical University in Prague, Technick\'a 4, CZ-16626 Prague, Czech Republic (e-mail: didier.henrion@laas.fr)}

\maketitle

\begin{abstract}                
HIFOO is a public-domain Matlab package initially designed for $H_{\infty}$ fixed-order controller synthesis,
using nonsmooth nonconvex optimization techniques. It was later on extended to multi-objective synthesis,
including strong and simultaneous stabilization under $H_{\infty}$ constraints. In this paper we describe a further extension
of HIFOO to $H_2$ performance criteria, making it possible to address mixed $H_2/H_{\infty}$ synthesis. We give
implementation details and report our extensive benchmark results.
\end{abstract}

\begin{center}
\end{center}


\section{Introduction}

HIFOO is a public-domain Matlab
package originally conceived during a stay of Michael Overton at the
Czech Technical University in Prague, Czech Republic, in the summer of 2005. HIFOO relies
upon HANSO, a general purpose implementation of an hybrid algorithm for
nonsmooth optimization, mixing standard quasi-Newton (BFGS) and gradient sampling techniques.
The acronym
HIFOO (pronounce [ha\i fu:]) stands for H-infinity Fixed-Order Optimization,
and the package is aimed at designing a stabilizing
linear controller of fixed-order for a linear plant in standard state-space configuration
while minimizing the $H_\infty$ norm of the closed-loop transfer function.

The first version of HIFOO was released and presented during the IFAC Symposium on
Robust Control Design in Toulouse, France in the summer of 2006, see \cite{hifoo2006},
based on the theoretical achievements reported in \cite{bhlo}.
HIFOO was later on extended to cope with multiple plant stabilization and multiple
conflicting objectives and the second major release of HIFOO was announced
during the IFAC Symposium on Robust Control Design in Haifa, Israel, in
the summer of 2009, see \cite{hifoo2009}.

Since then HIFOO has been used by various scholars and engineers. Benefiting
from feedback from users, we feel that it is now timely to extend HIFOO to
$H_2$ norm specifications. Indeed, $H_2$ optimal design, a generalization of
the well-known linear quadratic regulator design, is traditionally used in
modern control theory jointly with $H_{\infty}$ optimal design, see \cite{doyle}.
In particular, the versatile framework of mixed $H_2$/$H_\infty$ design
described e.g. in \cite{scherer}
is frequently used when designing high-performance control laws
for example in aerospace systems, see \cite{pirola}.
See also \cite{joris} for an application of the $H_2$ norm for smoothening
$H_{\infty}$ optimization.

The objective is this paper is to describe the extension of HIFOO to
$H_2$ norm specifications in such a way that users understand the basic mechanisms
underlying the package, and may be able to implement their own extensions
to fit their needs for their target applications. For example, the
algorithms of HIFOO can also be extended to cope with discrete-time
systems, pole placement specifications or time-delay systems.
On the HIFOO webpage
\begin{center}
\tt www.cs.nyu.edu/overton/software/hifoo
\end{center}
we are maintaining a list of publications reporting such extensions
and applications in engineering. The HIFOO and HANSO packages
can also be downloaded there.

\section{$H_2$ and $H_2/H_\infty$ synthesis}

\subsection{$H_2$ synthesis}
We use the standard state-space setup
\[
\begin{array}{rcl}
\dot{x} & = & Ax+B_1w+B_2u \\
z & = & C_1x+D_{11}w+D_{12}u \\
y & = & C_2x+D_{21}w+D_{22}u
\end{array}
\]
where $x$ contains the states, $u$ the physical (control) inputs, $y$ the
physical (measured) outputs, $w$ the performance inputs and $z$ the
performance outputs. Without loss of generality, we assume that
\[
D_{22} = 0
\]
otherwise we can use a linear change of variables on the system
inputs and outputs, see e.g. \cite{doyle}.

We want to design a controller
with state-space representation
\[
\begin{array}{rcl}
\dot{x}_K & = & A_K x_K + B_K y \\
u & = & C_K x_K + D_K y
\end{array}
\]
so that the closed-loop system equations become
\[
\begin{array}{rcl}
\dot{\bf x} & = & {\bf A}(k) {\bf x} + {\bf B}(k) w \\
z & = & {\bf C}(k) {\bf x} + {\bf D}(k) w
\end{array}
\]
in the extended state vector ${\bf x} = [x^T \:\: x_K^T]^T$ with matrices
\[
\begin{array}{rcl}
{\bf A}(k) & = & \left[\begin{array}{cc}
A+B_2D_KC_2 & B_2C_K \\ B_KC_2 & A_K
\end{array}\right] \\[1em]
{\bf B}(k) & = & \left[\begin{array}{c}
B_1+B_2D_KD_{21} \\ B_KD_{21}
\end{array}\right] \\[1em]
{\bf C}(k) & = & \left[\begin{array}{cc}
C_1+D_{12}D_KC_2 & B_KD_{21}
\end{array}\right] \\[.5em]
{\bf D}(k) & = & D_{11}+D_{12}D_KD_{21}
\end{array}
\]
depending affinely on the vector $k$ containing all
parameters in the controller matrices.

The $H_2$ norm of the closed-loop
transfer function $T(s)$ between input $w$ and output $z$
is finite only if matrix ${\bf A}$ is
asymptotically stable and if ${\bf D}$ is zero
(no direct feedthrough). This enforces the following
affine constraint on the $D_K$ controller matrix:
\begin{equation}\label{d22}
D_{11}+D_{12}D_KD_{21} = 0.
\end{equation}
We use the singular value decomposition to rewrite
this affine constraint in an explicit parametric
vector form, therefore reducing the number of
parameters in controller vector $k$.
If the above system of equations
has no solution, then there is no controller achieving
a finite $H_2$ norm.

In order to use the quasi-Newton optimization algorithms
of HANSO, we must provide a function evaluating the
$H_2$ norm in closed-loop and its gradient,
given controller parameters. Formulas can already be
found in the technical literature \cite{h2},
but they are reproduced here for
the reader's convenience. The (square of the) norm of the transfer
function $T(s)$ is given by
\[
f(k) = \|T(s)\|^2_2 = \mathrm{trace}\:({\bf C}X(k){\bf C}^T)
= \mathrm{trace}\:({\bf B}^TY(k){\bf B})
\]
where matrices $X(k)$ and $Y(k)$ solve the Lyapunov equations
\begin{equation}\label{lyap}
\begin{array}{rcl}
{\bf A}^T(k)X(k)+X(k){\bf A}(k)+{\bf C}^T(k){\bf C}(k) & = & 0, \\
{\bf A}(k)Y(k)+{\bf A}^T(k)Y(k)+{\bf B}(k){\bf B}^T(k) & = & 0
\end{array}
\end{equation}
and hence depend rationally on $K$.
The gradient of the $H_2$ norm with respect to controller
parameters $K$ is given by:
\[
\begin{array}{rcl}
\nabla_K f(k) & = & 2({\bf B}^T_2X(k)+{\bf D}^T_{12}{\bf C}(k))Y(k){\bf C}_2\\
& & + 2{\bf B}^T_2X(k){\bf B}(k){\bf D}^T_{21}
\end{array}
\]
upon defining the augmented system matrices
\[
\begin{array}{l}
{\bf B}_2 = \left[\begin{array}{cc}0&B_2\\1&0\end{array}\right], \quad
{\bf C}_2 = \left[\begin{array}{cc}0&1\\C_2&0\end{array}\right], \\[1em]
{\bf D}_{12} = [0 \:\: D_{12}], \quad {\bf D}_{21} = \left[\begin{array}{c}
0\\D_{21}\end{array}\right].
\end{array}
\]

\begin{figure}[h!]
\begin{center}
\includegraphics[width=8.4cm]{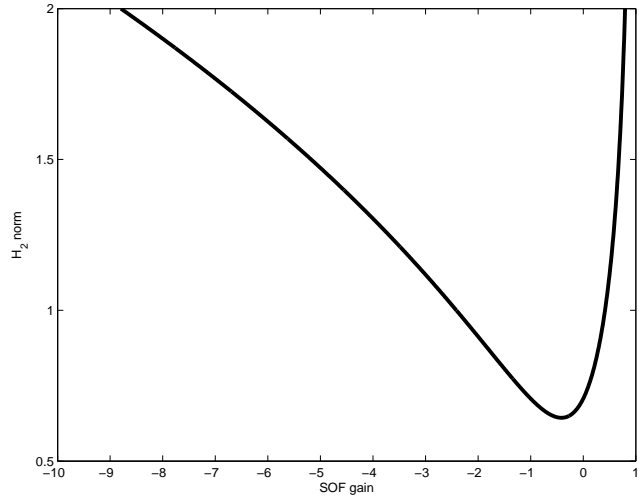}    
\caption{$H_2$ norm as a function of feedback gain.
\label{h2opt}}
\end{center}
\end{figure}

As an academic example for which the $H_2$ optimal controller can
be computed analytically, consider the system
\[
\begin{array}{rcl}
\dot{x} & = & -1+w+u \\
z & = & \left[\begin{array}{c}1\\0\end{array}\right] x +
\left[\begin{array}{c}0\\1\end{array}\right] u \\
y & = & x
\end{array}
\]
with a static controller
\[
u=ky
\]
with $k$ a real scalar to be found.
Closed-loop system matrices are
\[
{\bf A} = -1+k, \:\:
{\bf B} = 1, \:\:
{\bf C} = \left[\begin{array}{c}1\\k\end{array}\right].
\]
The first Lyapunov equation in (\ref{lyap}) reads
\[
2(-1+k)X(k)+1+k^2 = 0
\]
so the square of the $H_2$ norm is equal to
\[
f(k) = \frac{1+k^2}{2(1-k)}.
\]
For the gradient computation, we have to solve the second
Lyapunov equation in (\ref{lyap})
\[
2(-1+k)Y(k)+1 = 0
\]
and hence
\[
\nabla f(k) = \frac{1+2k-k^2}{2(1-k)^2}.
\]
This gradient vanishes at two points, one of which
violating the closed-loop stability condition $-1+k<0$.
The other point yields the optimal feedback gain
\[
k^* = 1-\sqrt{2} \approx -0.4142
\]
see Figure \ref{h2opt}. 

\begin{figure}[h!]
\begin{center}
\includegraphics[width=8.4cm]{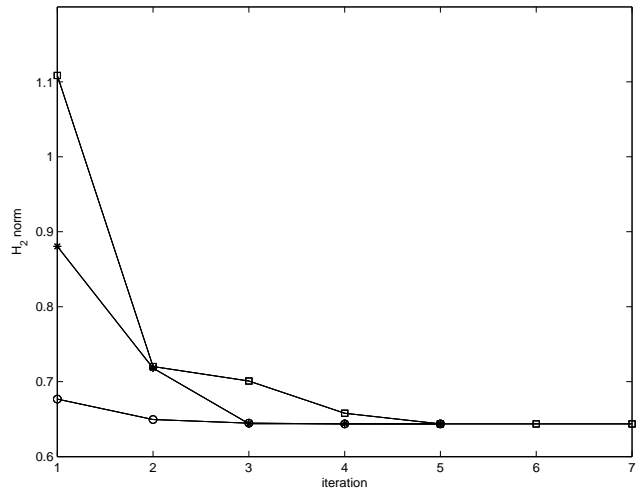}    
\caption{$H_2$ norm sequences optimized within HIFOO.
\label{h2its}}
\end{center}
\end{figure}

Using HIFOO with the input
sequence
\begin{verbatim}
P=struct('A',-1,'B1',1,'B2',1,...
'C1',[1;0],'C2',1,'D11',[0;0],...
'D12',[0;1],'D21',0,'D22',0);
options.prtlevel=2;
K=hifoo(P,'t',options)
\end{verbatim}
we generate the 3 sequences of optimized $H_2$ norms
displayed on Figure \ref{h2its}, yielding an
optimal $H_2$ norm of $0.6436$ consistent with
the analytic global minimum $\sqrt{\sqrt{2}-1}$.
Note the use of
the optional third input parameter specifying
a verbose printing level. Note also that the
sequences generated on your own computer may
differ since random starting points are used.

\subsection{Mixed $H_2/H_\infty$ synthesis}

\begin{figure}[h!]
\begin{center}
\includegraphics[width=6cm]{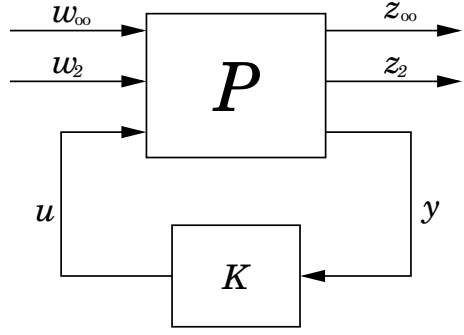}    
\caption{Standard feedback configuration for mixed $H_2/H_{\infty}$ synthesis.
\label{h2hinfgen}}
\end{center}
\end{figure}

One interesting feature of adding $H_2$ performance in the HIFOO package is the possibility to address the general mixed $H_2/H_\infty$ synthesis problem depicted on Figure \ref{h2hinfgen} where the open-loop plant is denoted by $P$ and the controller is denoted by $K$. 
A minimal state-space realization of the plant is given by
\[
P(s) := \left [\begin{array}{c|ccc}
A & B_{\infty} & B_{2} & B\\
\hline
C_{\infty} & D_{\infty} & {\bf 0} & D_{\infty u}\\
C_{2} & {\bf 0} & {\bf 0} & D_{2u}\\
C & D_{y\infty} & {\bf 0} & {\bf 0}
\end{array}\right ].
\]
The optimization problem
reads
\[
\begin{array}{ll}
\min_K & \|P_2(s)\|_2 \\
\mathrm{s.t.} & \|P_{\infty}(s)\|_{\infty} \leq \gamma_{\infty}
\end{array}
\]
where $P_2(s)$ is the transfer function between $H_2$ performance
signals $w_2$ and $z_2$, and $P_{\infty}(s)$ is the transfer function
between $H_{\infty}$ signals $w_{\infty}$ and $z_{\infty}$:
\begin{equation}
\begin{array}{lcl}
P_2(s) & := & \left [\begin{array}{c|cc}
A & B_2 & B\\
\hline
C_2 & {\bf 0} & D_{2u}\\
C & {\bf 0} & {\bf 0}
\end{array}\right ]\\ \\[-.5em]
P_\infty(s) & := & \left [\begin{array}{c|cc}
A & B_\infty & B\\
\hline
C_\infty & D_\infty & D_{\infty u}\\
C & D_{y\infty} & {\bf 0}
\end{array}\right ].
\end{array}
\label{plants-2-inf}
\end{equation}

An academic example for which the global optimal solution has been
calculated in \cite{cdc02} is used as an illustration for the mixed $H_2/H_\infty$ synthesis problem. Data for the model are given by
\[
\begin{array}{lll}
A=\left [\begin{array}{cc}
0 & 1\\
-1 & 0\end{array}\right ] & B=\left [\begin{array}{c}
0\\
1\end{array}\right ] & C=\left [\begin{array}{cc}
0 & 1\end{array}\right ]\\
  \\
C_{2}=\left [\begin{array}{cc}
1 & 0\\ 0 & 0\end{array}\right ] & B_{2}={\bf 1}_{2} & D_{2u}=\left [\begin{array}{c}
0\\
1\end{array}\right ]\\
  \\
C_{\infty}= \left [\begin{array}{cc}
0 & 1\end{array}\right ]& B_{\infty}= \left [\begin{array}{c}
1\\
 0\end{array}\right ] & D_{\infty u}=0\\
  \\
D_{\infty}=0 & D_{y\infty}=0 & D_{y2}={\bf 0}_{1\times 2}.\end{array}
\]
The analytical solution may be found by solving the following mathematical programming problem as in \cite{cdc02}
\begin{equation}
\begin{array}{lll}
& \displaystyle \min_{k} & \displaystyle J(k)\\
& \mathrm{s.t.} &\\
& & k<0\\
  & & f(k)\leq \gamma_\infty.
\end{array}
\label{pnl-h2hinf}
\end{equation}
For a non redundant mixed $H_2/H_\infty$ ($1<\gamma_\infty<\displaystyle\frac{3}{\sqrt{5}}$), the global optimal solution is 
\begin{equation}
\begin{array}{l}
k^{*}=-\sqrt{2-2\sqrt{1-1/\gamma^2}} \\
\\
 \|P_2\|_2=\alpha^{*}=\sqrt{\frac{4-3\sqrt{1-1/\gamma^{2}}}{\sqrt{2-2\sqrt{1-1/\gamma^2}}}}.
\end{array}
\label{solution-h2hinf}
\end{equation}

For $\gamma=1.2$, HIFOO gives the global optimal solution
\begin{equation}
\begin{array}{ll}
k^{*}=-0.9458\\ 
\|P_2\|_2=1.5735 & \|P_\infty\|_{\infty}=1.2
\end{array}
\label{gamma1.2}
\end{equation}
with the input sequence
\begin{verbatim}
P2=struct('A',[0 1;-1 0],'B1',eye(2),'B2',[0;1],...
'C1',[1 0;0 0],'C2',[0 1],'D11',zeros(2,2),...
'D12',[0;1],'D21',[0 0],'D22',0);

Pinf = struct('A',[0 1;-1 0],'B1',[1;0],...
'B2',[0;1],'C1',[0 1],'C2',[0 1],'D11',0,...
'D12',0,'D21',0,'D22',0);

K=hifoo({P2,Pinf},'th',[Inf,1.2]) 
\end{verbatim}

\section{Implementation details}

\subsection{Implementing $H_2$ norm into HIFOO}

In the main HIFOO function {\tt hifoomain.m}, we added an option {\tt 't'}
for $H_2$ norm specification, without affecting the existing features.
Proceeding this way, the $H_2$ norm can enter the objective function
or a performance constraint. The $H_2$ synthesis works in the same way
as the $H_\infty$ synthesis, using a stabilization phase followed by
an optimization phase.

A typical call of HIFOO for H2 static output feedback design is as follows:
\begin{verbatim}
K = hifoo('AC1','t')
\end{verbatim}
where {\tt AC1} refers to a problem of the COMPlib database, see
\cite{hifoo2006,hifoo2009}.

We added the function {\tt htwo.m} computing the $H_2$ norm and its gradient,
given the controller parameters. We had to pay special attention to the
linear system of equations arising from constraint (\ref{d22}).
When the user also specifies the controller structure, we have added
this constraint to the existing $H_2$ constraints. To do this we had
to change the way the controller structure was treated by HIFOO.

The above formulae for the computation of the $H_2$ norm and its gradient
are given for $D_{22} = 0$, i.e. zero feedthrough matrix. When this matrix is nonzero
we can use the same functions for synthesis, with some precautions. By considering
the shifted output $\tilde{y}=y-D_{22}u$, we recover the initial case with zero
feedthrough matrix. We compute the controller matrices $\hat{A}_K$, $\hat{B}_K$,
$\hat{C}_K$, $\hat{D}_K$ based
on the shifted output and then we apply a transformation on the controller matrix.
We obtain the final solution:
\begin{equation}
\begin{array}{rcl}
  A_K & = & \hat{A}_K - \hat{B}_K D_{22} (1 + \hat{D}_K D_{22})^{-1} \hat{C}_K \\
  B_K & = & \hat{B}_K (1 - D_{22} (1 + \hat{D}_K D_{22})^{-1} \hat{D}_K) \\
  C_K & = & (1 + \hat{D}_K D_{22})^{-1} \hat{C}_K \\
  D_K & = & (1 + \hat{D}_K D_{22})^{-1} \hat{D}_K.
  \end{array}
\end{equation}
Given the way this case is treated, multiple plant optimization only works
if the plants have the same feedthrough matrix. Note however that the case of
nonzero feedthrough matrix and imposed controller structure cannot
be treated by the current version of the program but could be the object of further development.

\subsection{Numerical linear algebra}

For $H_{\infty}$ norm optimization, HIFOO calls Matlab's
function {\tt eig} to check stability, returns {\tt inf}
if unstable, and otherwise calls the Control System Toolbox
function {\tt norm}, which proceeds by bisection on
successive computations of spectra of Hamiltonian matrices.
This latter function relies heavily
on system matrix scaling, on SLICOT routines, and its
is regularly updated and improved by The MathWorks Inc.
We observe experimentally that calling {\tt eig} once before
calling {\tt norm} is negligible (less than 5\%)
in terms of total computational cost.

For $H_2$ norm optimization, HIFOO calls {\tt eig} to check
stability, returns {\tt inf} if unstable, and otherwise
calls Matlab's {\tt lyap} function to
compute the norm and its gradient. Experimentally, we
observe that the time spent by {\tt eig} to check stability
is approximately 20\% of the time spent to solve the two
Lyapunov functions.

So a priori stability check is negligible for $H_{\infty}$
optimization and comparatively small but not negligible for $H_2$
optimization. In this latter case there is some room for improvement,
but since the overall objective of the HIFOO project is not
performance and speed but reliability, we decided to keep
the stability check for $H_2$ optimization.

\section{Benchmarking}

\subsection{$H_2$ synthesis}

We have extensively benchmarked HIFOO on problem instances
studied already in \cite{lmi09} with an LMI/randomized
algorithm. Since random starting points are used in HIFOO
we kept the best results over 10 attempts each with 3
starting points, with no computation time limit.
We ran the algorithms only on systems which
are not open-loop stable. For comparison, we took the best
results obtained in \cite{lmi09}.
In Table \ref{bench}, we use the following notations:
\begin{description}
\item $\star$: linear system (\ref{d22}) has a unique
solution which is not stabilizing
\item $\bullet$: linear system (\ref{d22}) has no solution
\item $+$: algorithm initialized with a stabilizing controller
\item $\dagger$: no stabilizing controller was found
\item $r$: rank assumptions on problem data are violated.
\end{description}
Also $n_x$, $n_u$, $n_y$ denote the number of states, inputs
and outputs. In addition to $H_2$ norms obtained by the LMI
algorithm and HIFOO, we also report for information
the $H_2$ norm achievable by full-order controller design
with HIFOO. Numerical values are reported to three
significant digits for space reasons.

In some cases (e.g. {\tt IH} and {\tt CSE2})
we observe that the norms achieved with a full-order
controller are greater than the norms achieved with a
static output feedback controller. This is due to the
difficulty of finding a good initial point in the
full-order case. A more practical approach, not pursued here,
consists in gradually increasing the order of the
controller, using the lower order controller found
at the previous step.

\begin{table}[!h]
\begin{center}\footnotesize
\caption{$H_2$ norm achieved for SOF controller design with LMI/randomized methods and HIFOO,
and full-order controller design with HIFOO.\label{bench}}
\begin{tabular}{@{}c|ccc|ccc@{}}
&$n_x$&$n_u$&$n_y$&SOF LMI&SOF HIFOO&full HIFOO\\ \hline
{\tt AC1} & 5 & 3 & 3 & 3.41e-7 & 1.46e-9 & 1.81e-15 \\
{\tt AC2} & 5 & 3 & 3 & 0.0503 & 0.0503 & 0.0491 \\
{\tt AC5} & 4 & 2 & 2 & 1470 & 1470 & 1340 \\
{\tt AC9} & 10 & 4 & 5 & r & 1.44 & 1.41 \\
{\tt AC10} & 55 & 2 & 2 & r & 27.8 ($+$) & $\dagger$ \\
{\tt AC11} & 5 & 2 & 4 & 3.94 & 3.94 & 3.64 \\
{\tt AC12} & 4 & 3 & 4 & r & 0.0202 & 5.00e-5 \\
{\tt AC13} & 28 & 3 & 4 & 132 & 132 & 106 \\
{\tt AC14} & 40 & 3 & 4 & r & $\star$ & 7.00 \\
{\tt AC18} & 10 & 2 & 2 & 19.7 & 19.7 & 18.6 \\
{\tt HE1} & 4 & 2 & 1 & 0.0954 & 0.0954 & 0.0857 \\
{\tt HE3} & 8 & 4 & 6 & r & 0.812 & 0.812 \\
{\tt HE4} & 8 & 4 & 6 & 21.7 & 20.8 & 18.6 \\
{\tt HE5} & 4 & 2 & 2 & r & $\star$ & 1.59 \\
{\tt HE6} & 20 & 4 & 6 & r & $\bullet$ & $\bullet$ \\
{\tt HE7} & 20 & 4 & 6 & r & $\bullet$ & $\bullet$ \\
{\tt DIS2} & 3 & 2 & 2 & 1.42 & 1.42 & 1.40 \\
{\tt DIS4} & 6 & 4 & 6 & 1.69 & 1.69 & 1.69 \\
{\tt DIS5} & 4 & 2 & 2 & r & $\star$ & 1280 \\
{\tt JE2} & 21 & 3 & 3 & 1010 & 961 & 623 \\
{\tt JE3} & 24 & 3 & 6 & r & $\bullet$ & $\bullet$ \\
{\tt REA1} & 4 & 2 & 3 & 1.82 & 1.82 & 1.50 \\
{\tt REA2} & 4 & 2 & 2 & 1.86 & 1.86 & 1.65 \\
{\tt REA3} & 12 & 1 & 3 & 12.1 & 12.1 & 9.91 \\
{\tt WEC1} & 10 & 3 & 4 & 7.36 & 7.36 & 5.69 \\
{\tt BDT2} & 82 & 4 & 4 & r & 0.795 & 0.655 \\
{\tt IH} & 21 & 11 & 10 & 1.66 & 1.54e-4 & 0.203 \\
{\tt CSE2} & 60 & 2 & 30 & 0.00890 & 0.00950 & 0.0133 \\
{\tt PAS} & 5 & 1 & 3 & 0.00920 & 0.00380 & 0.00197 \\
{\tt TF1} & 7 & 2 & 4 & r & 0.164 & 0.136 \\
{\tt TF2} & 7 & 2 & 3 & r & $\dagger$ & 10.9 \\
{\tt TF3} & 7 & 2 & 3 & r & 13.6 & 0.136 \\
{\tt NN1} & 3 & 1 & 2 & 41.8 & 41.8 & 35.0 \\
{\tt NN2} & 2 & 1 & 1 & 1.57 & 1.57 & 1.54 \\
{\tt NN5} & 7 & 1 & 2 & 142 & 142 & 82.4 \\
{\tt NN6} & 9 & 1 & 4 & 1350 & 1310 & 314 \\
{\tt NN7} & 9 & 1 & 4 & 133 & 133 & 84.2 \\
{\tt NN9} & 5 & 3 & 2 & r & 29.7 & 20.9 \\
{\tt NN12} & 6 & 2 & 2 & 18.9 & 18.9 & 10.9 \\
{\tt NN13} & 6 & 2 & 2 & r & $\bullet$ & $\bullet$ \\
{\tt NN14} & 6 & 2 & 2 & r & $\star$ & $\dagger$ \\
{\tt NN15} & 3 & 2 & 2 & 0.0485 & 0.0486 & 0.0480 \\
{\tt NN16} & 8 & 4 & 4 & 0.298 & 0.291 & 0.342 \\
{\tt NN17} & 3 & 2 & 1 & 9.46 & 9.46 & 3.87 \\
{\tt HF2D10} & 5 & 2 & 3 & 7.12e4 & 7.12e4 & 7.06e4 \\
{\tt HF2D11} & 5 & 2 & 3 & 8.51e4 & 8.51e4 & 8.51e4 \\
{\tt HF2D14} & 5 & 2 & 4 & 3.74e5 & 3.74e5 & 3.73e5 \\
{\tt HF2D15} & 5 & 2 & 4 & 2.97e5 & 2.97e5 & 2.84e5 \\
{\tt HF2D16} & 5 & 2 & 4 & 2.85e5 & 2.85e5 & 2.84e5 \\
{\tt HF2D17} & 5 & 2 & 4 & 3.76e5 & 3.76e5 & 3.75e5 \\
{\tt HF2D18} & 5 & 2 & 2 & 27.8 & 27.8 & 24.3 \\
{\tt TMD} & 6 & 2 & 4 & r & 1.36 & 1.32 \\
{\tt FS} & 5 & 1 & 3 & 1.69e4 & 1.69e4 & 1.83e4 \\
\end{tabular}
\end{center}
\end{table}

For the considered examples, HIFOO generally gives better results
than the randomized/LMI method of \cite{lmi09}. We also report
the performance achievable with a full-order controller designed
with HIFOO. We could not use the $H_2$ optimal synthesis functions of the
Control System Toolbox for Matlab as the technical assumptions (rank
conditions on systems data) under
which these functions are guaranteed to work are most of the time violated.

\subsection{$H_2$ synthesis for larger order systems}

For larger order systems, we compared our results with those of
\cite{aude} which are also based on nonsmooth optimization (labeled NSO).
In Table \ref{benchlarge} the column $n_k$ indicates the order
of the designed controller.
We observe that HIFOO yields better results, except for example {\tt CM4}
in the static output feedback case.

\begin{table}[h!]
\begin{center}
\caption{$H_2$ norm achieved with HIFOO compared with the
nonsmooth optimization method of \cite{aude}.\label{benchlarge}}
\begin{tabular}{c|ccc|c|cc}
&$n_x$&$n_u$&$n_y$&$n_k$&HIFOO&NSO\\\hline
{\tt AC14} & 40 & 4 & 3 & 1 & 21.4 & 21.4 \\
& & & & 10 & 7.00 & 8.10 \\
& & & & 20 & 7.00 & 7.56 \\
{\tt BDT2} & 82 & 4 & 4 & 0 & 0.791 & 0.794 \\
& & & & 10 & 0.598 & 0.789 \\
& & & & 41 & 0.585 & 0.779 \\
{\tt HF1} & 130 & 1 & 2 & 0 & 0.0582 & 0.0582 \\
& & & & 10 & 0.0581 & 0.0582 \\
& & & & 25 & 0.0581 & 0.0581 \\
{\tt CM4} & 240 & 1 & 2 & 0 & 61.0 & 0.926 \\
& & & & 50 & 0.933 & 0.938
\end{tabular}
\end{center}
\end{table}

\subsection{Mixed $H_2/H_{\infty}$ synthesis}

In this section we compare the results achieved with HIFOO
with those of \cite{aude} in the case of the mixed
$H_2/H_{\infty}$ synthesis problem, depicted on Figure
\ref{h2hinfgen}.

In Table \ref{benchmixed}
$n_k$ is the order of the controller and $\gamma_{\infty}$
is the level of $H_{\infty}$ performance (a constraint).
For problem dimensions
refer to Table \ref{benchlarge}.

We observe that HIFOO returns better or similar results than
the non-smooth optimization (NSO) method of \cite{aude},
except for problem {\tt CM4} in the static output feedback case.
Based on the CPU times of the NSO method gracefully provided
to us by Aude Rondepierre (not reported here), we must however
mention that HIFOO is typically much slower. This is not
surprising however since HIFOO is Matlab interpreted, contrary
to the NSO method which is compiled.

\begin{table}[h!]
\begin{center}
\caption{Mixed $H_2/H_{\infty}$ design with HIFOO compared with the
nonsmooth optimization method of \cite{aude}.\label{benchmixed}}
\begin{tabular}{c|cc|c@{\;}c|c@{\;}c}
& $n_k$ & $\gamma_{\infty}$ & $H_2$ HIFOO & $H_2$ NSO &
$H_{\infty}$ HIFOO & $H_{\infty}$ NSO \\ \hline
{\tt AC14} & 1 & 1000 & 21.4 & 21.4 & 230 & 231 \\
& 10 & 1000 & 7.01 & 8.78 & 100 & 101 \\
& 1 & 200 & 21.7 & 21.5 & 200 & 200 \\
& 20 & 200 & 7.08 & 7.99 & 100 & 100 \\
{\tt BDT2} & 0 & 10 & 0.790 & 0.804 & 0.908 & 1.06 \\
& 10 & 10 & 0.608 & 0.765 & 0.867 & 1.11 \\
& 0 & 0.8 & 0.919 & 0.791 & 0.943 & 0.800 \\
& 10 & 0.8 & 1.16 & 0.772 & 1.23 & 0.800 \\
& 41 & 0.8 & 1.24 & 0.789 & 2.32 & 0.800 \\
{\tt HF1} & 0 & 10 & 0.0582 & 0.0582 & 0.460 & 0.461 \\
& 0 & 0.45 & 0.0588 & 0.0588 & 0.450 & 0.450 \\
& 10 & 0.45 & 0.0586 & 0.0587 & 0.450 & 0.450 \\
& 25 & 0.45 & 0.0586 & 0.0587 & 0.450 & 0.450 \\
{\tt CM4} & 0 & 10 & 0.927 & 0.927 & 1.66 & 1.66 \\
& 0 & 1 & 0.986 & 0.984 & 1.00 & 1.00 \\
& 25 & 1 & 1.25 & 0.953 & 10.4 & 1.00
\end{tabular}
\end{center}
\end{table}

\section{Conclusion}

This paper documents the extension of HIFOO to $H_2$ performance. The resulting
new version 3.0 of HIFOO has been extensively benchmarked on $H_2$ and
$H_2/H_{\infty}$ minimization problems. We illustrated that HIFOO gives better
results than alternative methods for most of the considered benchmark problems.

HIFOO is an open-source public-domain
software that can be downloaded at
\begin{center}
\tt www.cs.nyu.edu/overton/software/hifoo
\end{center}
Feedback from users is welcome and significantly helps us improve the software
and our understanding of nonsmooth nonconvex optimization methods
applied to systems control.

Just before the completion of this work, Pierre Apkarian informed us that
several algorithms of nonsmooth optimization have now been implemented by The
MathWorks Inc. and will be released in the next version of the Robust
Control Toolbox for Matlab. Extensive comparison with HIFOO will
therefore be an interesting further research topic.

\section*{Acknowledgments}

This work benefited from feedback by Wim Michiels, Marc Millstone, Michael Overton and Aude Rondepierre.
The research of Suat Gumussoy was supported by the Belgian Programme on Interuniversity Poles
of Attraction, initiated by the Belgian State, Prime Minister's Office for Science, Technology
and Culture, and of the Optimization in Engineering Centre OPTEC.
The research of Didier Henrion was partly supported by project No.~103/10/0628 of the Grant Agency of the Czech Republic.
Denis Arzelier gratefully acknowledges Michael Overton and
the Courant Institute of Mathematical Sciences of NYU,
New York City, New York, USA, for hospitality, where this
work was initiated. Support for this work was provided in part
by the grant DMS-0714321 from the U.S. National Science
Foundation.


\begin{thebibliography}{99}

\bibitem{aude}
P. Apkarian, D. Noll and A. Rondepierre.
Mixed $H_2/H_{\infty}$ control via nonsmooth optimization.
SIAM Journal on Control and Optimization, 47(3):1516-1546, 2008.

\bibitem{pirola}
D. Arzelier, B. Cl\'ement, D. Peaucelle.
Multi-objective $H_2$/$H_\infty$/impulse-to-peak control of a space launch vehicle.
European Journal of Control, Vol. 12, No. 1, 2006.

\bibitem{cdc02}
D. Arzelier, D. Peaucelle. An iterative method for mixed H2/Hinfinity
synthesis via static output-feedback.
IEEE Conference on Decision and Control (CDC), Las Vegas, Nevada, 2002.

\bibitem{lmi09}
D. Arzelier, E.N. Gryazina, D. Peaucelle, B.T. Polyak. Mixed LMI/randomized
methods for static output feedback control design: stability and
Performance.  Technical Report LAAS-CNRS,  N°09640, September 2009.
Available at {\tt homepages.laas.fr/arzelier}

\bibitem{bhlo}
J. V. Burke, D. Henrion, A. S. Lewis and M. L. Overton.
Stabilization via nonsmooth, nonconvex optimization.
IEEE Transactions on Automatic Control, 51(11):1760-1769, 2006.

\bibitem{hifoo2006}
J. V. Burke, D. Henrion, A. S. Lewis and M. L. Overton. HIFOO - A Matlab Package for Fixed-order Controller Design and
H-infinity Optimization.  IFAC Symposium on Robust Control Design, Toulouse, France, 2006.

\bibitem{hifoo2009}
S. Gumussoy, D. Henrion, M. Millstone and M.L. Overton. Multiobjective Robust Control with HIFOO 2.0.
IFAC Symposium on Robust Control Design (ROCOND), Haifa, Israel, 2009.

\bibitem{h2}
T. Rautert and E. W. Sachs. Computational design of optimal output feedback
controllers, SIAM Journal on  Optimization, 7(3):837-852, 1997.

\bibitem{scherer}
C. W. Scherer. Multi-objective $H_2$/$H_\infty$ control.
IEEE Transactions on Automatic Control, 40:1054--1062, 1995.

\bibitem{joris}
J. Vanbiervliet, B. Vandereycken, W. Michiels, S. Vandewalle and M. Diehl.
Smoothed spectral abscissa for robust stability optimization.
SIAM Journal on Optimization, 20(1):156-171, 2009.

\bibitem{doyle}
K. Zhou, J. Doyle and K. Glover.
Robust and Optimal Control,
Prentice-Hall, 1996.

\end{thebibliography}
\end{document}